\documentclass[11pt]{article}
\usepackage{amsmath, amssymb, cite, colordvi}
\usepackage{graphicx, ifpdf}

\newtheorem{thm}{Theorem}[section]
\newtheorem{cor}[thm]{Corollary}
\makeatletter \@addtoreset{equation}{section}
\makeatletter \@addtoreset{figure}{section}

\def\qed{\hfill \rule{4pt}{7pt}}
\def\pf{\noindent {\it Proof.}\quad}

\makeatletter
\def\ExtendSymbol#1#2#3#4#5{\ext@arrow 0099{\arrowfill@#1#2#3}{#4}{#5}}
\makeatother

\title{\bf\LARGE The Method of Combinatorial Telescoping}
\author{William Y.C. Chen\raisebox{5pt}{\scriptsize 1},
Qing-Hu Hou\raisebox{5pt}{\scriptsize 2} and Lisa H. Sun\raisebox{5pt}{\scriptsize 3}}
\date{Center for Combinatorics, LPMC-TJKLC\\
 Nankai University\\
Tianjin 300071, P.R. China \\
\vspace{15pt}
\raisebox{5pt}{\scriptsize 1\,}chen@nankai.edu.cn,
\raisebox{5pt}{\scriptsize 2\,}hou@nankai.edu.cn, \raisebox{5pt}{\scriptsize 3\,}sunhui@nankai.edu.cn}

\begin{document}
\maketitle

\noindent {\bf Abstract.}
We present a method for proving $q$-series identities by  combinatorial
telescoping, in the sense that one can transform a bijection or a classification
of combinatorial objects into a telescoping relation. We shall illustrate this method
by giving a combinatorial proof of Watson's identity which implies
the Rogers-Ramanujan identities.



\section{Introduction}

The main objective of this paper is to present the method of combinatorial
telescoping for proving $q$-series identities. A benchmark of this approach is
the classical identity of Watson which implies the Rogers-Ramanujan
identities.

There have been many combinatorial proofs of the Rogers-Ramanujan identities.
Schur \cite{Pak} provided an involution for the following  identity
which is equivalent to the first Rogers-Ramanujan identity:
\[
\prod_{k=1}^\infty (1-q^k) \left(1 + \sum_{k=1}^\infty \frac{q^{k^2}}{(1-q)(1-q^2) \cdots (1-q^k)} \right)  = \sum_{k=-\infty}^{\infty} (-1)^k q^{k(5k-1)/2}.
\]
Andrews \cite{Andrews} proved the Rogers-Ramanujan identities by introducing the
notion of $k$-partitions. Garsia and Milne \cite{Garsia} gave a  bijection by using the involution principle. Bressoud and Zeilberger \cite{BZ1,BZ2} provided a different involution principle proof based on an algebraic proof due to Bressoud \cite{Bre}.
Boulet and Pak \cite{Boulet} found a combinatorial proof
 which relies on the symmetry properties of a generalization of Dyson's rank.

Let us consider a summation of the following form
\begin{equation}\label{eq-f}
\sum\limits_{k=0}^\infty (-1)^k f(k).
\end{equation}
Suppose that $f(k)$ is a weighted count
of a set $A_k$, that is,
\[
f(k)=\sum\limits_{\alpha \in A_k} w(\alpha).
\]
Motivated by the idea of
  creative telescoping of Zeilberger \cite{Zeil}, we aim  to
find sets $B_k$ and $H_k$ with a weight assignment $w$  such that there is a weight preserving bijection
\begin{equation} \label{ab}
\phi_k \colon A_k \longrightarrow B_k \cup H_k \cup H_{k+1},
\end{equation}
where $\cup$ stands for disjoint union.
Since  $\phi_k$ and $\phi_{k+1}$
 are weight preserving, both $\phi_k^{-1}(H_{k+1})$ and $\phi_{k+1}^{-1}(H_{k+1})$ have the same  weight as $H_{k+1}$.  Realizing that $\phi_k^{-1}(H_{k+1}) \subseteq A_k$ and $\phi_{k+1}^{-1}(H_{k+1}) \subseteq A_{k+1}$, they cancel  each other in the sum \eqref{eq-f}.  More precisely, if we set
 \[ g(k) =\sum_{\alpha\in B_k} w(\alpha) \quad \mbox{and} \quad
      h(k)=\sum_{\alpha \in H_k} w(\alpha),\]
      then the bijection (\ref{ab})
 implies  that
 \begin{equation}
 f(k)=g(k)+h(k)+h(k+1).
 \end{equation}
 To see that the above equation is indeed a telescoping relation with respect to the sum (\ref{eq-f}),
  let
\[
 f'(k)=(-1)^k f(k), \quad g'(k)= (-1)^k g(k), \quad h'(k)=(-1)^k h(k).\]
Thus we have
\begin{equation}\label{teles}
  f'(k) =  g'(k) + h'(k) - h'(k+1).
\end{equation}
Just like the conditions for  the creative telescoping, we suppose that $H_0=\emptyset$ and $H_k$ vanishes for sufficiently large $k$. Summing   \eqref{teles} over $k$, we deduce the following relation
\begin{equation}\label{eqAtoB}
\sum_{k=0}^\infty (-1)^k f(k)=\sum_{k=0}^\infty (-1)^k g(k),
\end{equation}
which is often an identity we wish to establish.

The above approach to proving an identity like
(\ref{eqAtoB}) is called combinatorial telescoping.
It can be seen that the bijections $\phi_k$ lead to a
correspondence between $A=\mathop{\bigcup}\limits_{k=0}^\infty A_k$ and $B=\mathop{\bigcup}\limits_{k=0}^\infty B_k$ after the cancelations of $H_k$'s.
To be more specific,  we can derive a bijection
\[
\phi \colon A \setminus \bigcup_{k=0}^\infty \phi_k^{-1}(H_k \cup H_{k+1}) \longrightarrow  B
\]
and an involution
\[
\psi \colon \bigcup_{k=0}^\infty \phi_k^{-1}(H_k \cup H_{k+1}) \longrightarrow \bigcup_{k=0}^\infty \phi_k^{-1}(H_k \cup H_{k+1}),
\]
given by $\phi(\alpha) = \phi_k(\alpha)$ if $\alpha \in A_k$ and
\[
\psi(\alpha) =   \begin{cases}
                  \phi_{k-1}^{-1} \phi_k(\alpha), \quad \mbox{if}\ \alpha \in \phi_k^{-1}(H_k), \\[5pt]
                  \phi_{k+1}^{-1} \phi_k(\alpha), \quad \mbox{if}\ \alpha \in \phi_k^{-1}(H_{k+1}).
                  \end{cases}
\]

In the examples of this paper, the set $A_k$ is of the following form
\[
A_k=\bigcup_{n=0}^\infty A_{n,k}.
\]
Fix an integer $n$, for any nonnegative integer $k$, we can establish a bijection
$\phi_{n,k} $ such that the corresponding set $B_{n,k}$ is related to $A_{n,k}, A_{n-1,k}, \ldots, A_{n-r,k}$ for an  integer $r$.
Let
\[
F_{n,k} = \sum_{\alpha\in A_{n,k}} w(\alpha)
\]
be a weighted count of the set $A_{n,k}$, and let
\[
F_n = \sum_{k=0}^\infty (-1)^k F_{n,k}.
\]
By \eqref{eqAtoB}, the bijections $\{\phi_{n,k}\}_{k=0}^\infty$ imply a recurrence relation of $F_n$, which leads to an explicit expression $u(n)$ for $F_n$ by iteration. Finally, we deduce the following identity
\begin{equation}\label{i2}
\sum_{k=0}^\infty (-1)^k f(k) = \sum_{k=0}^\infty  (-1)^k \sum_{n=0}^\infty F_{n,k} = \sum_{n=0}^\infty F_n  = \sum_{n=0}^\infty u(n).
\end{equation}

As a simple example, one can easily give a combinatorial telescoping proof of
the classical identity of Gauss, see also, \cite{chl03,Pak04,lee10}:
\[
\sum_{k=0}^n (-1)^k {n \brack k} = \begin{cases}
0, & \mbox{$n$ odd,} \\
(1-q)(1-q^3) \cdots (1-q^{n-1}), & \mbox{$n$ even.}
\end{cases}
\]
Let us consider the following reformulation
\begin{equation}\label{eq-Gauss}
\sum_{k=0}^n (-1)^k \frac{1}{(q;q)_k (q;q)_{n-k}} = \begin{cases}
0, & \mbox{$n$ odd,} \\[5pt]
\displaystyle \frac{1}{(1-q^2)(1-q^4) \cdots (1-q^n)}, & \mbox{$n$ even.}
\end{cases}
\end{equation}
Let
\[
P_{n,k} = \{(\lambda, \mu) \colon \lambda_1 \le k, \ \mu_1 \le n-k\},
\]
where $\lambda$ and $\mu$ are partitions, and let
\[
H_{n,k} = \{(\lambda, \mu) \in P_{n,k} \colon m_{k}(\lambda) < m_{n-k}(\mu) \},
\]
where $m_k(\lambda)$ denotes the number of occurrences of the part $k$ in $\lambda$ and we
adopt the convention that $m_0(\lambda)=+\infty$. By definition, $H_{n,k}=\emptyset$ for $k=0$ or $k>n$.
For any integers $n \ge 1$ and $k \ge 0$, we shall construct
a bijection
\[
\phi_{n,k} \colon P_{n,k} \longrightarrow \{0,n,2n,\ldots\} \times P_{n-2,k} \cup H_{n,k} \cup H_{n,k+1}.
\]

Let $(\lambda,\mu) \in P_{n,k}$. If $m_k(\lambda) < m_{n-k}(\mu)$, then $(\lambda, \mu) \in H_{n,k}$. In this case, $\phi_{n,k}\big((\lambda,\mu)\big)=(\lambda,\mu)$. If $m_k(\lambda) \ge m_{n-k}(\mu)$, we let $m_{n-k}(\mu) =t$. In this case, if $\mu_{t+1} = n-1-k$, we increase each of the first $t$ parts of $\lambda$ by one and decrease each of the first $t$ parts of $\mu$ by one. It is easily seen that the resulting pair of partitions $(\lambda', \mu')$ belongs to $H_{n,k+1}$ and we set  $\phi_{n,k} \big((\lambda, \mu)\big)=(\lambda', \mu')$. Finally, if $\mu_{t+1} \le n-2-k$, then we set
\[
\phi_{n,k} \big((\lambda, \mu)\big) = \big(tn, ( \hat{\lambda}, \hat{\mu}) \big) \in \{0,n,2n,\ldots\} \times P_{n-2,k},
\]
where $
\hat{\lambda} = (\lambda_{t+1}, \lambda_{t+2}, \ldots)$ and  $\hat{\mu}=(\mu_{t+1}, \mu_{t+2}, \ldots)
$
are the partitions obtained from $\lambda$ and $\mu$ by removing the first $t$ parts.
Define the weight function $w$ on $P_{n,k}$ and $\{0,n,2n,\ldots\} \times P_{n-2,k}$ as follows
\[
w(\lambda, \mu) = q^{|\lambda|+|\mu|}, \quad \mbox{and} \quad w(tn, (\lambda, \mu)) = q^{tn+|\lambda|+|\mu|},
\]
where $|\lambda|=\lambda_1+\lambda_2+\cdots$.
It can be checked that $\phi_{n,k}$ is weight
preserving. Hence we obtain the following recurrence relation
\begin{equation}\label{fn}
F_n(q) = \frac{1}{1-q^{n}} F_{n-2}(q),
\end{equation}
where $F_n(q)$ denotes the sum on the left hand side of \eqref{eq-Gauss}. By iteration
of (\ref{fn}), we arrive at \eqref{eq-Gauss}.

It should be noted that the bijections $\phi_{n,k}$ lead to an involution on $P_{n,k}$,
which can be considered as a variation of the
involution given by Chen, Hou and Lascoux \cite{chl03}.

In Section~\ref{sec-Schur}, we use the idea of combinatorial telescoping to give a proof of Watson's identity \cite{Watson} in the following form, see also \cite[Section
2.7]{Gasper},
\begin{equation}\label{eq-Watson}
\sum_{k=0}^\infty (-1)^k \frac{1-aq^{2k}}{(q;q)_k (aq^k;q)_\infty}a^{2k}q^{k(5k-1)/2}= \sum_{n=0}^\infty \frac{a^nq^{n^2}}{(q;q)_n},
\end{equation}
where
\[
(a;q)_k=(1-a)(1-aq) \cdots (1-aq^{k-1}), \quad \mbox{and} \quad (a;q)_\infty=\prod_{i=0}^\infty (1-a q^i).
\]
Setting $a=1$, Watson's identity reduces to Schur's identity \cite{Boulet}
\[
\frac{1}{(q;q)_\infty} \sum_{k=-\infty}^\infty (-1)^k q^{k(5k-1)/2} = \sum_{n=0}^\infty \frac{q^{n^2}}{(q;q)_n}.
\]
Applying Jacobi's triple product identity to the left hand side, we are led to the first Rogers-Ramanujan identity. Similarly, setting $a=q$ in Watson's identity yields the second Rogers-Ramanujan identity.

Here is a sketch of the proof.
Assume that
the $k$-th summand regardless of the sign on the left hand side of \eqref{eq-Watson}
is the weight of a set $P_k$.
      We further divide $P_k$ into a disjoint union of subsets $P_{n,k}, n=0,1,\ldots$,
      by considering the expansion of the summand in the parameter $a$.
      For a positive integer $n$ and a nonnegative integer $k$, we can construct a bijection
\begin{equation}\label{rec}
\phi_{n,k} \colon P_{n,k} \to \{n\} \times P_{n,k} \cup \{2n-1\} \times P_{n-1,k}  \cup H_{n,k} \cup H_{n,k+1}.
\end{equation}
Let
\[
F_n(a,q) = \sum_{k=0}^\infty (-1)^k \sum_{\alpha \in P_{n,k}} w(\alpha).
\]
The  bijections $\phi_{n,k}$ yield a recurrence relation
\[
F_n(a,q) = q^n F_n(a,q) + a q^{2n-1} F_{n-1}(a,q), \quad \; n \ge 1.
\]
By iteration, we find
 that $F_n(a,q)= a^n q^{n^2}/(q;q)_n$, and hence \eqref{eq-Watson} holds.

As another example, it can be seen that the method of combinatorial telescoping also applies to Sylvester's identity \cite{Sylv}
\begin{equation}\label{eq-Sylv}
\sum_{k=0}^\infty (-1)^k q^{k(3k+1)/2} x^k \frac{1-xq^{2k+1}}{(q;q)_k (xq^{k+1};q)_\infty} = 1.
\end{equation}
This identity has been investigated by Andrews \cite{Andrews76, Andrews}.

\section{ Watson's identity}\label{sec-Schur}

In this section, we shall use Watson's identity as an example to illustrate
the idea of combinatorial telescoping.
Let us recall some definitions concerning partitions.
A {\it partition} is a non-increasing finite sequence of positive integers $\lambda=(\lambda_1,\ldots,\lambda_\ell)$. The integers $\lambda_i$ are called the {\it parts} of $\lambda$.
The sum of parts and the number of parts are denoted by $|\lambda|=\lambda_1+\cdots + \lambda_\ell$ and $\ell(\lambda)=\ell$, respectively.  The number of $k$-parts in $\lambda$ is denoted by $m_{k}(\lambda)$. The special partition with no parts is denoted by $\varnothing$. We shall use diagrams to represent
partitions and use  columns instead of rows to represent parts.

Set \begin{equation}
P_k=\{(\tau, \lambda, \mu)  \colon \tau=(k^{2k},k-1,\ldots, 2, 1), \ \lambda_{\ell(\lambda)} \ge k, \, \lambda_i \not= 2k, \ \mu_1 \le k
\},
\end{equation}
where $k^{2k}$ denotes $2k$ occurrences of a part $k$.
In other words, $\tau$ is a trapezoid partition with $|\tau|=k(5k-1)/2$, $\lambda$ is a partition with parts at least $k$ but not equal to $2k$, and $\mu$ is a partition with parts at most $k$. In particular, we have $P_0=\{(\varnothing, \lambda, \varnothing)\}$.
It is clear that  the $k$-th summand of the left hand side of (\ref{eq-Watson}) without sign
  can be viewed as the
 weight of $P_k$, that is,
\[
\sum_{(\tau,\, \lambda, \,\mu) \in P_k} a^{\ell(\lambda)+2k} q^{|\tau| + |\lambda| + |\mu|}.
\]

\begin{figure}[ht]
\centering
\setlength{\unitlength}{0.05 mm}%
  \begin{picture}(1365.9, 547.8)(0,0)
  \put(0,0){\includegraphics{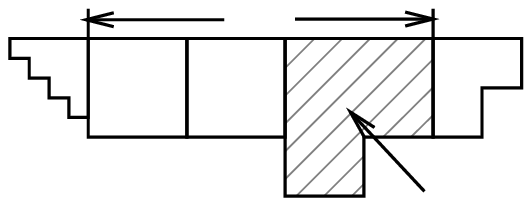}}
  \put(533.42,446.00){\fontsize{8.53}{10.24}\selectfont \makebox(30.0, 60.0)[l]{$n$\strut}}
  \put(220.48,314.11){\fontsize{8.53}{10.24}\selectfont \makebox(30.0, 60.0)[l]{$k$\strut}}
  \put(749.18,302.10){\fontsize{8.53}{10.24}\selectfont \makebox(60.0, 60.0)[l]{$\lambda$\strut}}
  \put(100.56,147.93){\fontsize{11.38}{13.66}\selectfont \makebox(200.0, 80.0)[l]{$\underbrace{\qquad\qquad\qquad}$\strut}}
  \put(302.62,38.20){\fontsize{8.53}{10.24}\selectfont \makebox(90.0, 60.0)[l]{$\tau$\strut}}
  \put(853.67,41.56){\fontsize{8.53}{10.24}\selectfont \makebox(472.2, 60.0)[l]{$\ge k$ but without $2k$\strut}}
  \put(942.41,310.75){\fontsize{8.53}{10.24}\selectfont \makebox(60.0, 60.0)[l]{$\mu$\strut}}
  \end{picture}%
\caption{\label{pic1}The diagram $(\tau,\lambda,\mu) \in P_{n,k}$}
\end{figure}

According to the exponent of $a$ in the above definition, we divide $P_k$ into a disjoint union of subsets \begin{equation}\label{defP}
P_{n,k} = \{ (\tau,\lambda,\mu) \in P_k \colon \ell(\lambda) = n-2k \},
\end{equation}
with $P_{n,0}=\{(\varnothing,\lambda,\varnothing) \in P_0 \colon \ell(\lambda) = n\}$ and $P_{n,k}=\emptyset$ for $n<2k$.
The elements of $P_{n,k}$ are illustrated in Figure~\ref{pic1}.

We have the following combinatorial telescoping relation for  $P_{n,k}$.

\begin{thm}\label{mainlem}
Let
\begin{equation}\label{defI}
H_{n,k}=\{(\tau,\lambda,\mu) \in P_{n,k} \colon m_{k}(\lambda)+2 > m_{k}(\mu) \}.
\end{equation}
Then, for any positive integer $n$ and  any nonnegative integer $k$, there is a
bijection
\begin{equation}\label{P-ct}
\phi_{n,k} \colon P_{n,k} \longrightarrow \{n\} \times P_{n,k} \cup \{2n-1\} \times P_{n-1,k} \cup H_{n,k}
\cup H_{n,k+1}.
\end{equation}
\end{thm}

\pf The bijection is essentially a classification of $P_{n,k}$ into four cases.
Let $(\tau,\lambda,\mu)$  be a $3$-tuple of
 partitions in $ P_{n,k}$.

\noindent
Case 1. $m_{k}(\lambda) +2> m_{k}(\mu)$. In this case,  $(\tau,\lambda,\mu) \in H_{n,k}$ and the image of $(\tau,\lambda,\mu)$ is defined to be itself.

\noindent
Case 2. $m_{k}(\lambda)+2 \leq m_{k}(\mu)$ and $m_{2k+1}(\lambda)=0$.
Denote the set of $3$-tuples $(\tau, \lambda, \mu)$
in this case by $U_{n,k}$. Note that
\[
U_{n,0}=\{(\varnothing,\lambda,\varnothing) \in P_{n,0} \colon  m_1(\lambda)=0\}.
\]
Since $m_{k}(\mu) \ge m_{k}(\lambda)+2$, we can remove $(m_{k}(\lambda)+2)$ $k$-parts from $\mu$ to generate a partition $\mu'$. In the meantime, we change each $k$-part of $\lambda$ into a $2k$-part in order to obtain a partition $\lambda'$ whose minimal part is strictly greater than $k$.

\begin{figure}[ht]
\centering
  \setlength{\unitlength}{0.05 mm}%
  \begin{picture}(1120.0, 504.4)(0,0)
  \put(0,0){\includegraphics{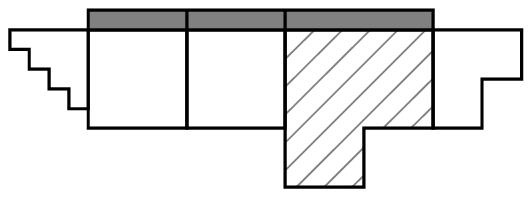}}
  \put(220.48,242.62){\fontsize{8.53}{10.24}\selectfont \makebox(30.0, 60.0)[l]{$k$\strut}}
  \put(749.18,230.61){\fontsize{8.53}{10.24}\selectfont \makebox(120.0, 60.0)[l]{$\lambda''$\strut}}
  \put(942.41,239.25){\fontsize{8.53}{10.24}\selectfont \makebox(90.0, 60.0)[l]{$\mu'$\strut}}
  \put(280.48,402.62){\fontsize{8.53}{10.24}\selectfont \makebox(30.0, 60.0)[l]{$k$\strut}}
  \put(480.48,402.62){\fontsize{8.53}{10.24}\selectfont \makebox(30.0, 60.0)[l]{$k$\strut}}
  \put(660.48,402.62){\fontsize{8.53}{10.24}\selectfont \makebox(120.0, 60.0)[l]{$n-2k$\strut}}
  \end{picture}%
\caption{\label{piccase2}The resulting partition under the bijection $\varphi_1$.}
\end{figure}

 Next, we decrease  each part of $\lambda'$ by one in order to produce
 a partition $\lambda''$ whose minimal part is greater than or equal to $k$. Since $\lambda$ contains no parts equal to $2k+1$, we see that
 $\lambda''$ contains no parts equal to $2k$. Thus we  obtain a bijection
$\varphi_1\colon  U_{n,k}   \rightarrow    \{n\} \times P_{n,k}$ defined by
                  $ (\tau, \lambda, \mu)   \mapsto
                  (n, (\tau, \lambda'', \mu'))$.
                  This case is illustrated by Figure~\ref{piccase2}.

\begin{figure}[ht]
\centering
  \setlength{\unitlength}{0.05 mm}%
  \begin{picture}(1120.0, 544.4)(0,0)
  \put(0,0){\includegraphics{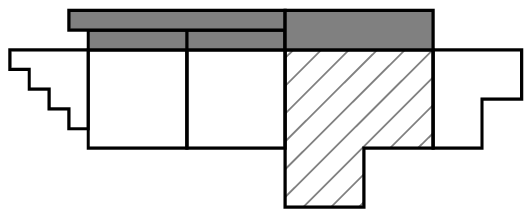}}
  \put(220.48,242.62){\fontsize{8.53}{10.24}\selectfont \makebox(30.0, 60.0)[l]{$k$\strut}}
  \put(749.18,230.61){\fontsize{8.53}{10.24}\selectfont \makebox(120.0, 60.0)[l]{$\lambda''$\strut}}
  \put(942.41,239.25){\fontsize{8.53}{10.24}\selectfont \makebox(90.0, 60.0)[l]{$\mu'$\strut}}
  \put(620.48,442.62){\fontsize{8.53}{10.24}\selectfont \makebox(180.0, 60.0)[l]{$n-2k-1$\strut}}
  \put(300.48,442.62){\fontsize{8.53}{10.24}\selectfont \makebox(120.0, 60.0)[l]{$2k+1$\strut}}
  \end{picture}%
\caption{\label{piccase3}The resulting partition under the bijection $\varphi_2$.}
\end{figure}

\noindent
Case 3. $m_{k}(\lambda)+2 \leq m_{k}(\mu)$, $m_{2k+1}(\lambda)>0$ and $m_{k+1}(\lambda)+ m_{2k+2}(\lambda)=0$. Denote the set of $3$-tuples $(\tau, \lambda, \mu)$
in this case by $V_{n,k}$. We remark that when $k=0$, one $1$-part is regarded as a $(2k+1)$-part and the other $1$-parts are regarded as $(k+1)$-parts so that
\[
V_{n,0} = \{(\varnothing, \lambda, \varnothing) \in P_{n,0} \colon \mbox{$m_1(\lambda)=1$ and $m_2(\lambda)=0$} \}.
\]
Let $\lambda',\mu'$ be given as in Case 2. We can remove one  $(2k+1)$-part from $\lambda'$ and decrease  each of the remaining parts by two in order to
 obtain $\lambda''$.
This leads to a bijection $\varphi_2\colon  V_{n,k}  \rightarrow   \{2n-1\} \times P_{n-1,k}$
as given by $(\tau, \lambda, \mu)   \mapsto
                  (2n-1, (\tau, \lambda'', \mu'))$.
 See Figure~\ref{piccase3} for  an illustration.

\begin{figure}[ht]
\centering
  \setlength{\unitlength}{0.05 mm}%
  \begin{picture}(1493.1, 547.8)(0,0)
  \put(0,0){\includegraphics{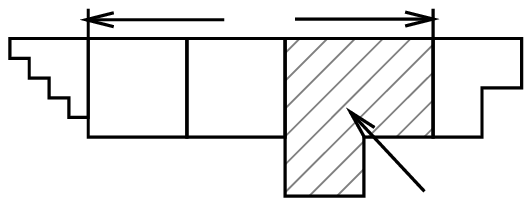}}
  \put(533.42,446.00){\fontsize{8.53}{10.24}\selectfont \makebox(30.0, 60.0)[l]{$n$\strut}}
  \put(220.48,314.11){\fontsize{8.53}{10.24}\selectfont \makebox(90.0, 60.0)[l]{$k+1$\strut}}
  \put(749.18,302.10){\fontsize{8.53}{10.24}\selectfont \makebox(150.0, 60.0)[l]{$\lambda'''$\strut}}
  \put(100.56,147.93){\fontsize{11.38}{13.66}\selectfont \makebox(200.0, 80.0)[l]{$\underbrace{\qquad\qquad\qquad}$\strut}}
  \put(302.62,38.20){\fontsize{8.53}{10.24}\selectfont \makebox(120.0, 60.0)[l]{$\tau'$\strut}}
  \put(853.67,41.56){\fontsize{8.53}{10.24}\selectfont \makebox(599.4, 60.0)[l]{$\ge k+1$ and without $2k+2$\strut}}
  \put(942.41,310.75){\fontsize{8.53}{10.24}\selectfont \makebox(120.0, 60.0)[l]{$\mu''$\strut}}
  \end{picture}%
\caption{\label{piccase4}The resulting partition under the bijection $\varphi_3$.}
\end{figure}

\noindent
Case 4. $m_{k}(\lambda)+2 \leq m_{k}(\mu)$, $m_{2k+1}(\lambda)>0$ and $m_{k+1}(\lambda)+ m_{2k+2}(\lambda)>0$. Denote the set of $3$-tuples $(\tau, \lambda, \mu)$
in this case by $W_{n,k}$. As in Case 3, we have
\[
W_{n,0}=\{(\varnothing,\lambda,\varnothing) \in P_{n,0} \colon \mbox{$m_1(\lambda)>0$ and $m_1(\lambda)+m_2(\lambda)>1$} \}.
\]
Let $\lambda',\mu'$ be given as in Case 2. We can change each $(2k+2)$-part of $\lambda'$ to a $(k+1)$-part and add $m_{2k+2}(\lambda')$ $(k+1)$-parts to $\mu'$. Denote the resulting partitions by $\lambda''$ and $\mu''$. Then we have
\begin{equation}\label{case-2.2}
m_{k+1}(\lambda'') =m_{k+1}(\lambda)+ m_{2k+2}(\lambda)>0, \quad
m_{k+1}(\mu'')=m_{2k+2}(\lambda).
\end{equation}
Now remove one $(k+1)$-part and one $(2k+1)$-part from $\lambda''$ to obtain $\lambda'''$. By \eqref{case-2.2}, we find
\[
m_{k+1}(\lambda''') = m_{k+1}(\lambda'')-1 \ge m_{k+1}(\mu'') - 1.
\]
Moreover, it is clear that
 \[|\lambda|+|\mu| = 2k + (k+1) + (2k+1) + |\lambda'''|+|\mu''|.\] Let $\tau'$ be the trapezoid partition of size $k+1$. So we obtain a bijection
$  \varphi_3\colon  W_{n,k} \rightarrow   H_{n,k+1}$ defined by
                  $(\tau, \lambda, \mu)  \mapsto
                  (\tau', \lambda''', \mu'')$.
This case is illustrated in Figure~\ref{piccase4}. \qed

Assign a weight function $w$ on $P_{n,k}, \{n\} \times P_{n,k}$ and $\{2n-1\} \times P_{n-1,k}$ as follows:
\begin{align*}
&w \big(\tau,\lambda,\mu\big) =a^n q^{|\tau|+|\lambda|+|\mu|}, \\
&w\big( n, (\tau,\lambda,\mu) \big)=q^n \cdot a^n  q^{|\tau|+|\lambda|+|\mu|},\\
&w\big( 2n-1, (\tau,\lambda,\mu) \big)=aq^{2n-1}\cdot a^{n-1} q^{|\tau|+|\lambda|+|\mu|}.
\end{align*}
Observe that the bijections $\varphi_1$, $\varphi_2$ and $\varphi_3$ are weight preserving.
In addition, $H_{n,0}=\emptyset$ and $H_{n,k}=\emptyset$ for $k>\frac{n}{2}$. Thus the bijections $\phi_{n,k}$ immediately lead to a recurrence relation of $F_n(a,q)$ defined as follows.

\begin{cor}\label{rec-Fn}
Let
\begin{equation}\label{def-Fn}
F_n(a,q) = \sum_{k=0}^\infty (-1)^k \sum_{(\tau,\lambda,\mu) \in P_{n,k}} a^n q^{|\tau| + |\lambda| + |\mu|}.
\end{equation}
Then, for any positive integer $n$, we have
\begin{equation}
\label{rec-F}
F_n(a,q) = q^n F_n(a,q) + a q^{2n-1} F_{n-1}(a,q).
\end{equation}
\end{cor}
Since $F_0(a,q) = 1$, by iteration we find that
\[
F_n(a,q) = \frac{aq^{2n-1}}{1-q^n} F_{n-1}(a,q)
= \frac{a^2 q^{4n-4}}{(1-q^n)(1-q^{n-1})} F_{n-2}(a,q)
= \cdots = \frac{a^n q^{n^2}}{(q;q)_n}.
\]
Summing over $n$, we arrive at Watson's identity \eqref{eq-Watson}.

\section{Sylvester's identity}\label{remark}

\begin{figure}[ht]
\centering
  \setlength{\unitlength}{0.05 mm}%
  \begin{picture}(980.0, 556.3)(0,0)
  \put(0,0){\includegraphics{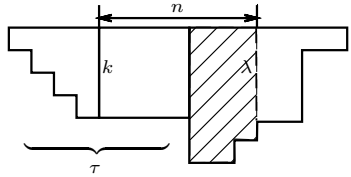}}
  \put(470.47,466.83){\fontsize{6.83}{8.19}\selectfont \makebox(24.0, 48.0)[l]{$n$\strut}}
  \put(287.22,305.03){\fontsize{6.83}{8.19}\selectfont \makebox(24.0, 48.0)[l]{$k$\strut}}
  \put(655.19,306.90){\fontsize{6.83}{8.19}\selectfont \makebox(48.0, 48.0)[l]{$\lambda$\strut}}
  \put(254.16,38.56){\fontsize{6.83}{8.19}\selectfont \makebox(72.0, 48.0)[l]{$\tau$\strut}}
  \put(81.34,113.51){\fontsize{6.83}{8.19}\selectfont \makebox(144.0, 48.0)[l]{$\underbrace{\rule{55pt}{0pt}}$\strut}}
  \end{picture}%
\caption{\label{Sylv}The diagram of $(\tau,\lambda) \in Q_{n,k}$.}
\end{figure}

In this section, we describe the approach of combinatorial telescoping for Sylvester's identity \eqref{eq-Sylv}.
Define
\[
Q_{n,k} = \{(\tau, \lambda) \colon \tau=(k^{k+1}, k-1, \ldots, 2, 1), \lambda_i \not= 2k+1, m_{>k}(\lambda) = n-k \},
\]
where $m_{> k}(\lambda)$ denotes the number of parts of $\lambda$ which are greater than $k$. See Figure~\ref{Sylv} for an illustration. In particular, we have
\[
Q_{n,0}= \{(\varnothing, \lambda) \colon  \lambda_i \not= 1, \ell(\lambda) = n\}.
\]

Let
\[
H_{n,k} = \{ (\tau,\lambda) \in Q_{n,k} \colon m_{k+1}(\lambda) \ge m_k (\lambda) \}.
\]
Then, for each positive integer $n$ and each nonnegative integer $k$, we have a bijection
\[
\phi_{n,k}: Q_{n,k} \longrightarrow \{n\}\times Q_{n,k} \cup H_{n,k} \cup H_{n,k+1},
\]
which is a classification of $Q_{n,k}$ into three cases. Let $(\tau,\lambda) \in Q_{n,k}$.

\noindent
Case 1. $m_{k+1}(\lambda) \ge m_k(\lambda)$. In this case,  $(\tau, \lambda) \in H_{n,k}$ and the image of $(\tau,\lambda)$ under $\phi_{n,k}$ is defined to be itself.

\noindent
Case 2. $m_{k+1}(\lambda) < m_k(\lambda)$ and $m_{2k+2}(\lambda)=0$. Denote the set of pairs $(\tau, \lambda)$ in this case by $U_{n,k}$. We remove one $k$-part from $\lambda$. Then, for each  $(k+1)$-part of $\lambda$, we can add it to a $k$-part to form a $(2k+1)$-part. Finally, we decrease each part greater than $k+1$ by one to generate a partition $\lambda'$.  Since $m_{2k+2}(\lambda)=0$, we see that  $(\tau, \lambda')\in Q_{n,k}$. So we obtain a bijection $\varphi_1 \colon U_{n,k} \to \{n\}\times Q_{n,k}$ given by $(\tau, \lambda) \mapsto(n,(\tau, \lambda'))$.

\noindent
Case 3. $m_{k+1}(\lambda) < m_k(\lambda)$ and $m_{2k+2}(\lambda)>0$. Denote the set of   pairs $(\tau, \lambda)$ in this case by $V_{n,k}$. We first remove one $k$-part and one $(2k+2)$-part from $\lambda$ and add them to $\tau$ to form a partition $\tau'$. Here $\tau'$ is a trapezoid partition of size $k+1$. Then for each  $(k+1)$-part of $\lambda$ we combine it with a $k$-part to form a $(2k+1)$-part. Finally we decompose each $(2k+3)$-part of $\lambda$ into a $(k+1)$-part and a $(k+2)$-part to form a partition $\lambda'$. Since $m_{2k+3}(\lambda')=0$, we obtain a bijection $\varphi_2 \colon V_{n,k} \to H_{n,k+1}$ defined by $(\tau, \lambda) \mapsto (\tau', \lambda')$.

It is not difficult to see that Sylvester's identity follows from the bijections $\phi_{n,k}$. Let
\[
I_n(q) = \sum_{k=0}^\infty (-1)^k \sum_{(\tau,\lambda) \in Q_{n,k}} q^{|\tau|+|\lambda|}.
\]
Noting that $H_{n,0}=\emptyset$ because of the definition $m_0(\lambda)=+\infty$, the bijections $\phi_{n,k}$ lead to the recurrence relation \[ I_n(q) = q^n I_n(q),\]
which implies that $I_n(q)=0$ for $n \ge 1$. Clearly $I_0(q)=1$, and hence Sylvester's identity holds.

To conclude this paper, we notice that both Watson's identity and Sylvester's identity can be verified by employing  the $q$-Zeilberger algorithm for infinite $q$-series
 developed by
Chen, Hou and Mu \cite{Chen-08}. Let
\[
f(a) = \sum_{k=0}^\infty (-1)^k \frac{(1-aq^{2k})}{(q;q)_k (aq^k;q)_\infty} a^{2k} q^{k(5k-1)/2}.
\]
Denote the $k$-th summand of $f(a)$ by $F_k(a)$. The $q$-Zeilberger algorithm  gives that
\begin{equation}\label{rec-fa}
F_k(a)-F_k(aq)-aqF_k(aq^2)=H_{k+1}(a)-H_k(a),
\end{equation}
where
\[
H_k(a)=(-1)^k\frac{(-1-q^k+aq^{2k})}{(q;q)_{k-1}(aq^k;q)_\infty} a^{2k} q^{k(5k-1)/2}.
\]
Summing (\ref{rec-fa}) over $k$, we find that \[ f(a) = f(aq) + aqf(aq^2).\] Extracting the coefficients of  $a^n$ leads to the same recurrence relation as \eqref{rec-F}.
It is easily checked that the right hand side of \eqref{eq-Watson} satisfies the same recursion. By Theorem~3.1 of Chen, Hou and Mu \cite{Chen-08},
one sees that
\eqref{eq-Watson} holds for any $a$ provided that it is valid for the trivial case $a=0$. Similarly, let
\[
f(x) = \sum_{k=0}^\infty (-1)^k q^{k(3k+1)/2} x^k \frac{1-xq^{2k+1}}{(q;q)_k (xq^{k+1};q)_\infty}.
\]
The $q$-Zeilberger algorithm gives that
\begin{equation}\label{fg}
F_k(x)-F_k(xq)=H_{k+1}(x)-H_k(x),
\end{equation}
where $F_k(x)$ is the $k$-th summand of $f(x)$ and
\begin{equation}\label{hk}
H_k(x)=(-1)^{k+1} \frac{q^{k(3k+1)/2} x^k}{(q;q)_{k-1} (xq^{k+1};q)_\infty}.
\end{equation}
Summing   (\ref{fg}) over $k$, we deduce that
$f(x)=f(xq)$, which implies $f(x)=1$.

\vskip 15pt \noindent {\small {\bf Acknowledgments.}  We wish to thank the referees for helpful comments. This
work was supported by the National Science Foundation, the PCSIRT project,  the Project NCET-09-0479, and the Fundamental Research Funds for  Universities of the Ministry of Education of China.}

\end{document}